\documentclass[12pt]{article}

\usepackage{amssymb}
\usepackage{amsthm}
\usepackage{amsmath}
\usepackage{amsfonts}
\usepackage{hyperref}
\usepackage{color}
\usepackage[T1]{fontenc}
\title{On certain variant of strongly nonlinear interpolation inequality in
dimension $n$ }

\author{Tomasz Choczewski$^{\rm a}$
\ and Agnieszka Ka\l{}amajska$^{\rm a,b,}$\thanks{Corresponding author. Email: A.Kalamajska@mimuw.edu.pl
\vspace{6pt}}\hspace{5pt}$^{,}$\footnote{{Supported by Warsaw Center of Mathematics and Computer Science (KNOW) at Institute of Mathematics,
 Polish Academy of Sciences, http://wcmcs.edu.pl.} }\\
\vspace{-6pt}
\small{$^{a}${\em{Faculty of Mathematics, Informatics, and Mechanics, University of Warsaw, Warszawa, Poland}}};\\
\small{$^{\rm b}${\em{Institute of Mathematics Polish Academy of Sciences, Warszawa, Poland}}}}
\renewcommand{\it}{\sl}
\renewcommand{\em}{\sl}

plus 6pt minus 12pt
\newcommand{\barint}{
         \rule[.036in]{.12in}{.009in}\kern-.16in
          \displaystyle\int  }
\def\r{{\mathbb{R}}}
\def\n{{\mathbb{N}}}
\def\rn{{\mathbb{R}^{n}}}

\begin{document}

\newtheorem{theo}{\bf Theorem}[section]
\newtheorem{coro}{\bf Corollary}[section]
\newtheorem{lem}{\bf Lemma}[section]
\newtheorem{rem}{\bf Remark}[section]
\newtheorem{defi}{\bf Definition}[section]
\newtheorem{ex}{\bf Example}[section]
\newtheorem{fact}{\bf Fact}[section]
\newtheorem{prop}{\bf Proposition}[section]
\newtheorem{prob}{\bf Problem}[section]

\makeatletter \@addtoreset{equation}{section}
\renewcommand{\theequation}{\thesection.\arabic{equation}}
\makeatother

\newcommand{\ds}{\displaystyle}
\newcommand{\ts}{\textstyle}
\newcommand{\ol}{\overline}
\newcommand{\wt}{\widetilde}
\newcommand{\ck}{{\cal K}}
\newcommand{\ve}{\varepsilon}
\newcommand{\vp}{\varphi}
\newcommand{\pa}{\partial}
\newcommand{\rp}{\mathbb{R}_+}
\newcommand{\hh}{\tilde{h}}
\newcommand{\HH}{\tilde{H}}
\newcommand{\ct}{{\cal T}}

\maketitle

%\bigskip {\it Dedicated to the memory of Professor Marek Burnat}

\smallskip
  {\small Key words and phrases: interpolation inequalities, multiplicative inequalities, Sobolev spaces\\

  }

  \begin{abstract}
  \noindent
We obtain the inequality $$\int_{\Omega}|\nabla u(x)|^ph(u(x))dx\leq C(n,p)\int_{\Omega} \left( \sqrt{ |\nabla^{(2)} u(x)||{\cal T}_{h,C}(u(x))|}\right)^{p}h(u(x))dx,$$ where
$\Omega\subseteq \rn$ and $n\ge 2$,  $u:\Omega\rightarrow \r$ is  in certain subset in second order Sobolev space $W^{2,1}_{loc}(\Omega)$,  $\nabla^{(2)} u$ is the Hessian matrix of $u$,
${\cal T}_{h,C}(u)$ is certain transformation of the continuous function $h(\cdot)$. Such inequality is the generalization of similar inequality holding in one dimension, obtained earlier by second author and Peszek.
\end{abstract}

\section{Introduction}
The purpose of this paper is to obtain the $n$-dimensional variant of the following inequality:
\begin{eqnarray}\label{jedwym}
\int_{(a,b)}|u^{'}(x)|^ph(u(x))dx\leq C_p\int_{(a,b)}\sqrt{\left|u^{''}(x){\cal T}_h(u(x))\right|}^{p}h(u(x))dx,
\end{eqnarray}
where  $(a,b)$ is an interval, $u\in W_{loc}^{2,1}((a,b))$  and  obeys some additional assumptions, $h(\cdot )$ the a given continuous  function,
${\cal T}_h(\cdot )$ is a certain  transform of $h(\cdot )$.  When for example $h(t) = t^\alpha, \alpha >-1,$ then ${\cal T}_{h}(\lambda)=\frac{1}{\alpha +1}\lambda$, so it is proportional to $\lambda$, see Remark \ref{przyprzyk}, where we deal with ${\cal T}_{h,0}(\lambda)$.
The constant $C_p$ does not depend on $u$.

The above inequality was obtained by second author and Peszek in \cite{akjp}.
It implies the inequality:
\[
\left( \int_\r |f^{'}|^{p}h(f) dx\right)^{\frac{2}{p}} \le (p-1) \left(  \int_\r |{\cal T}_h(f)|^qh(f)dx\right)^{\frac{1}{q}} \left( \int_\r |f^{''}|^rh(f)dx\right)^{\frac{1}{r}},
\]
where $q\ge \frac{p}{2}$, $\frac{2}{p}=\frac{1}{q}+\frac{1}{r}$. After the substitution of $h\equiv 1$ we obtain  the classical Gagliardo-Nirenberg inequality (\cite{ga,n1}):
\[
\left( \int_\r |f^{'}|^{p} dx\right)^{\frac{2}{p}} \le C \left(  \int_\r |f|^qdx\right)^{\frac{1}{q}} \left( \int_\r |f^{''}|^rdx\right)^{\frac{1}{r}},
\] where in our case one has to assume that $q\ge \frac{p}{2}$. However, inequality \eqref{jedwym} within $q\ge \frac{p}{2}$ is more general.

Inequalities of the type \eqref{jedwym}, were generalized later to various settings  in \cite{cfk,km,kp2}. They were applied to results about regularity and asymptotic behavior in nonlinear eigenvalue problems for singular pde's like for example  the Emden-Fowler equation
\cite{km,akjp}:
\(
u^{''}(x)= g(x)u(x)^\alpha ,\  \alpha\in {\bf {R}}, \ g(x)\in L^p((a,b))
\), which appear for example in electricity theory, fluid dynamics or mathematical biology. The example equation with $g(t)=  t^{\frac{1}{2}}$, $\alpha = {\frac{3}{2}}$ appears in
in the model found independently in 1927 by Thomas and Fermi  to determine the electrical potential in an isolated neutral atom (\cite{fermi}, \cite{tomas}).
Inequalities like  \eqref{jedwym} were also applied to
 the study of regularity of solutions of the   Cucker-Smale equation with  singular weight, \cite{pesz}.
In the more theoretical approach, they were used to obtain  certain  generalization of isoperimetric inequalities and capacitary estimates due to Mazya  \cite{kp2} in the Orlicz setting.

 All results obtained in the papers \cite{km,akjp,kp2}
were dealing only with the case of inequalities for one-variable function and
 $p\ge
2$. In the most recent paper \cite{cfk},  the authors obtained the generalization of the one dimensional inequality
\eqref{jedwym} to the case $1<p<2$. Unexpectedly, it appears that in the new inequality instead of the operator
${\cal T}_h(u(x))$ involves certain nonlocal  operator, i.e. such a one which uses all values of $u$ in the interval $(a,b)$, not just at $x$.

We obtained the inequalities in the form:
\begin{equation}\label{cozrobilismy}
\int_{\Omega}|\nabla u(x)|^ph(u(x))dx\leq C(n,p)\int_{\Omega} \left( \sqrt{ |\nabla^{(2)} u(x)||{\cal T}_{h,C}(u(x))|}\right)^{p}h(u(x))dx,
\end{equation}
where  $C(n,p)=\left( p-1+\sqrt{n-1} \right)^{\frac{p}{2}}$, $\Omega\subseteq \rn$ and $n\ge 2$,  $u:\Omega\rightarrow \r$ belongs to certain subset in the Sobolev space $W^{2,1}_{loc}(\Omega)$,  $\nabla^{(2)} u$ is the Hessian matrix of $u$ and
${\cal T}_{h,C}(u)$ is certain transformation of the continuous function $h(\cdot)$ (see Definition \ref{hh}), under certain additional assumptions (see Theorems \ref{main1} and \ref{main2}).

Let us indicate three earlier sources, where
 the variants of such inequality could be found. Namely, the inspiration for the authors of \cite{akjp} to derive the inequality \eqref{jedwym},
were the estimates due to Mazja (\cite{ma}, Section 8.2.1):
\[
\int_{{\rm supp} f^{'}} \left(\frac{ |f^{'}|}{f^{\frac{1}{2}}}\right)^pdx \le \left( \frac{p-1}{|1-\frac{1}{2} p|}\right)^{\frac{p}{2}}
 \int_{\bf R} |f^{''}|^{{p}}dx,
\]
where ${\rm supp} f^{'}$ is the support of $f^{'}$, which are the special case of \eqref{jedwym} when one considers $h(s)=s^{-p/2}$. They
were applied in  \cite{ma} to  the second order izoperimetric inequalities and capacitary estimates in second order Sobolev classes.

Opial obtained inequalities:
\begin{equation*}\label{opial11}
\int_a^b |yy^{'}|dx\le K \left( \int_a^b |y^{''}|^p dx\right)^{\frac{2}{p}},
\end{equation*}
known as second order Opial inequalities, holding on a compact interval $[a,b]$,
where
\[
y\in BC_0=\{ y\in W^{2,p}((a,b)): y(a)=y(b)=y^{'}(b)=0\},
\]
see e.g. see \cite{opial} and  e.g. \cite{bloom,bbch} for further related issues.

\noindent
In another source \cite{kp1}, the authors obtained the inequality:
\begin{eqnarray*}
\int_\rn G(|\nabla u|)dx&\le& C\int_\rn G(|u||\nabla^{(2)}u|)dx,
 \end{eqnarray*}
dealing with convex function $G$.

Let us mention that the passage from  the inequality  \eqref{jedwym} to the inequality \eqref{cozrobilismy} is not so direct and it requires quite delicate arguments. For example now we apply integrals over the
boundaries of Lipshitz domains and  the coarea formulae.

We hope that the new derived inequalities can be applied  the regularity theory for the nonlinear eigenvalue problems in singular elliptic pde's in the similar way as it was done for functions of one variable.

\section{Preliminaries and notation}

\subsection{Basic notation}
In the sequel we assume that $\Omega\subseteq \r^n$ is an open domain, $n\in \n$. We use standard definition of strongly Lipschitz domain i.e. of the domain of class ${\cal C}^{0,1}$, see e.g. \cite{ma}, Section 1.1.9.
 We use the standard notation: $C_0^\infty (\Omega)$ to denote smooth functions with compact support, $W^{m,p}(\Omega)$ and $W^{m,p}_{loc}(\Omega)$
to denote the  global and local Sobolev functions defined on $\Omega$, respectively, while $C_0^\infty (\Omega, \r^k)$, $W^{m,p}(\Omega, \r^k)$  and $W_{loc}^{m,p}(\Omega,\r^k)$ will denote their vectorial counterparts.
%If $I\subseteq \mathbf{R}$ is an interval, by $AC(I)$ we denote functions which are absolutely continuous on $I$.
By $d\sigma$ we denote the $n-1$-dimensional Hausdorff measure.
%If $C\subseteq \r$ and $f$ is defined on $C$
%by $f\chi_C$ we denote an extension of $f$ by zero outside set $C$.
For $p>1$ we set
\begin{eqnarray}\label{fip}
\Phi_p(\lambda)=\left\{
\begin{array}{ccc}
|\lambda|^{p-2}\lambda & {\rm if}& \lambda\neq 0\\
0&{\rm if} & \lambda =0,
\end{array}
\right.
\ \ \  \Phi_p:  \mathbf{R}^n\rightarrow \mathbf{R}^n,
\end{eqnarray}
and we write $p^{'}:=p/(p-1)$.
When $u:\in W^{1,1}_{loc}(\Omega,\r^m)$,
by $Du$ we denote the matrix of its differential  $\left(\frac{\partial u^i (x)}{\partial x_j}\right)_{i\in \{ 1,\dots,n\}, j\in \{ 1,\dots,m\} }$. When $u\in W^{2,1}_{loc}(\Omega)$,
      by $\nabla^{(2)} u(x)$   we denote the Hessian of
 $u$, i.e. the matrix $\left(\frac{\partial^{2}u (x)}{\partial x_i\partial x_j}\right)_{i,j\in \{ 1,\dots,n\} }$.
If $A$ is a vector or matrix, by $|A|$ we denote its Euclidean norm. The same notation may be used to denote Lebesgue's measure of the measurable subset in $\rn$.

\subsection{Properties of absolutely continuous functions}

 We will need the following variant of Nikodym ACL
Characterization Theorem, which can be found
e.g. in \cite{ma}, Section 1.1.3.

\begin{theo}[Nikodym ACL Characterization Theorem]\label{nikos}
 \begin{description}
\item[i)] Let  $u\in~W^{1,1}_{\rm loc}(\rn)$. Then for every $i\in \{
1,\dots ,n\}$ and for almost every $a\in \r^{i-1}\times \{
0\}\times \r^{n-i}$ the function
\begin{equation}\label{tratt}
\r\ni t\mapsto u(a+te_i)
\end{equation}
is locally absolutely continuous on $\r$. In particular for almost
every point $x\in\r^n$ the distributional derivative $D_iu(x)$ is
the same as the usual derivative at  $x$.
\item[ii)]
Assume that $u\in L^1_{loc}(\Omega)$ and
for every $i\in \{ 1,\dots,n\}$ and for almost every $a\in
\r^{n-i-1}\times \{ 0\}\times \r^i$ the function in (\ref{tratt}) is
locally absolutely continuous on $\r$ and all the derivatives
$D_iu$ computed almost everywhere
are locally integrable on $\r^n$. Then $u$ belongs to
$W^{1,1}_{loc}(\r^n)$.
\item[iii)]
Let $1\le p\le \infty$,   $u\in W^{1,1}_{loc}(\rn)$ and $\Omega\subset \r^n$ be an open
subset. Then $u$ belongs to $W^{1,p}(\Omega)$ if and only if $u\in L^p(\Omega)$ and
every derivative  $D_iu$ computed almost everywhere,  belongs to the space $L^{p}(\Omega)$.
\end{description}
\end{theo}

\noindent

The following lemma can be easily proven by using the very definition of absolutely continuous functions.
\begin{lem}\label{lcomp}
\begin{description}
\item[i)]
If $f: [-R,R]\rightarrow \r$ is absolutely continuous with values in the interval $[\alpha,\beta]$ and $L: [\alpha,\beta] \rightarrow\r$ is a  Lipschitz function, then
the function $(L\circ f)(x):=L(f(x))$ is absolutely continuous on $[-R,R]$;
\item[ii)]
If $L: [-R,R]\rightarrow \r$ is a Lipschitz function with values in the interval $[\alpha,\beta]$ and $f: [\alpha,\beta] \rightarrow\r$ is  absolutely continuous, then
the function $(f\circ L)(x):=f(L(x))$ is absolutely continuous on $[-R,R]$.
\end{description}
\end{lem}

\subsection{Transformations of the nonlinear weight}
\noindent
We will use the following definition, which is the special variant of definition introduced in  \cite{akjp}.

\begin{defi}[transforms of the nonlinear weight $h$]\label{hh}
Let $0< B\le \infty$, $h:(0,B)\rightarrow (0,\infty)$ be a continuous  function which is  integrable on  $(0,\lambda)$ for every $\lambda <B$, $C\in\r$,
and let $H_C:[0,B)\rightarrow \mathbf{R}$ be the locally absolutely continuous primitive of $h$ extended to $0$, given by
\[
H_C(\lambda):= \int_0^\lambda h(s)ds-C,\ \ \lambda \in [0,B).
\]
We define
the transform of $h$, ${\cal T}_{h,C}: (0,B)\rightarrow (0,\infty)$ by
\[
{\cal T}_{h,C}(\lambda):=
\frac{H_C(\lambda)}{h(\lambda)}, \ \ \lambda \in (0,B).
\]
\end{defi}
In the case when  $C=0$ we omit it from the notation, i.e. we write $H_0=:H,{\cal T}_{h,0} =: {\cal T}_{h}$.
Note that $h$ and  ${\cal T}_{h,C}$ might not be defined at $0$ or $B$ in case $B<\infty$.

\begin{rem}\label{przyprzyk}\rm
Simple examples of the admitted weights $h(\cdot)$ and their transformations can be found among power weights. Namely, when we chose $h(\lambda)=\lambda^\theta$ where $\theta >-1$, then clearly ${\cal T}_{h,0}(\lambda)=(1+\theta)^{-1}\lambda$, so it is proportional to the identity function. Similarly, ${\cal T}_{h,C}(\lambda)$  can be estimated by the proportional to the identity function  when   for example we consider $C=0$ and we can deduce that (*): ${\cal T}_{h,0}(\lambda)= \frac{H_0(\lambda)}{h(\lambda)}\le C\lambda$, with some general constant $C$. This is always the case when for $H_0$ is convex, so when $h$ is nondecreasing.
The example of power function shows that even in the case of decreasing functions $h(\cdot)$, the estimate (*) can still hold true.
\end{rem}

\section{Presentation of main results}%\label{three}

\noindent
Our goal is to prove the following $n$-dimensional variant of inequality  \eqref{jedwym}.

\begin{theo}\label{main1}
Let $n\ge 2$, $\Omega\subseteq \mathbf{R}^n$  be the bounded domain,
$2\le p<\infty$, $0<B\le \infty$, $h:(0,B)\rightarrow (0,\infty)$ and $H_C, {\cal T}_{h,C}$  be as in Definition \ref{hh}. Moreover, let the following assumptions be satisfied:
\begin{description}
\item[$(u)$:] $u\in  W^{2,p/2}_{loc}(\Omega)\cap C({\Omega})$
 and $0<u<B$ in $\Omega$;
\item[$(u,\Omega)$:] There exists sequence $\{ \Omega_k\}_{k\in \mathbf{N}}$  of bounded  subdomains of $\Omega$ of class ${\cal C}^{0,1}$ such that $\bar{\Omega}_k\subseteq \Omega$, $\bigcup_k\Omega_k=\Omega$ and
\begin{eqnarray}\label{aproksym}
\lim_{k\to\infty}\int_{\partial\Omega_k}\Phi_p(\nabla u(x))\cdot n(x) H_C(u(x))d\sigma (x)\in [-\infty,0],
\end{eqnarray}
where $n(x)$ denotes unit outer normal vector to $\partial\Omega_k$, defined for $\sigma$ almost all $x\in \partial\Omega_k$.
\end{description}
Then we have
\begin{eqnarray}\label{niermain1}
\int_{\Omega}|\nabla u(x)|^ph(u(x))dx\leq C(n,p)\int_{\Omega} \left( \sqrt{ |\nabla^{(2)} u(x)||{\cal T}_{h,C}(u(x))|}\right)^{p}h(u(x))dx,
\end{eqnarray}
where  $C(n,p)=\left( p-1+\sqrt{n-1} \right)^{\frac{p}{2}}$.
\end{theo}

Our next statement is the special variant of Theorem \ref{main1}, when we assume that $\Omega$ is of class ${\cal C}^{0,1}$ and we use some additional assumptions on $u$.

\begin{theo}\label{main2}
Let $n\ge 2$, $\Omega\subseteq \mathbf{R}^n$  be the bounded Lipschitz domain,
$2\le p<\infty$, $0<B\le \infty$, $h:(0,B)\rightarrow (0,\infty)$, $H_C, {\cal T}_{h,C}$  be as in Definition \ref{hh},
 $u\in  W^{2,p/2}(\Omega)\cap C(\bar{\Omega})$,
  $0<u<B$ a.e. in $\Omega$ and
  \begin{eqnarray}\label{aproksym1}
\int_{\partial\Omega}\Phi_p(\nabla u(x))\cdot n(x) H_C(u(x))d\sigma (x)\in [-\infty,0],
\end{eqnarray}
where $n(x)$ denotes unit outer normal vector to $\partial\Omega$, defined for $\sigma$ almost all $x\in \partial\Omega$.

Then
\begin{eqnarray}\label{niermain111}
\int_{\Omega}|\nabla u(x)|^ph(u(x))dx\leq C(n,p)\int_{\Omega} \left( \sqrt{ |\nabla^{(2)} u(x)||{\cal T}_{h,C}(u(x))|}\right)^{p}h(u(x))dx,
\end{eqnarray}
where  $C(n,p)=\left( p-1+\sqrt{n-1} \right)^{\frac{p}{2}}$.
\end{theo}

Interval $(0,B)$ in the assumptions of Theorems \ref{main1} and \ref{main2} can be changed to
any interval $(A,B)$ where $-\infty< A<B\le \infty$. Our most general result reads as follows.

\begin{theo}\label{main3}\rm
Let $n\ge 2$, $\Omega\subseteq \mathbf{R}^n$  be the bounded domain,
$2\le p<\infty$, $-\infty < A<B\le \infty$, $h: (A,B)\rightarrow (0,\infty)$ is continuous and integrable near $A$, $ C\in \mathbf{R}$,
\[
H_{C,A}(\lambda):=\int_A^{\lambda+A} h(\tau)d\tau-C,\ \ \ \
{\cal T}_{h,C,A}(\lambda) :=\frac{\int_A^{\lambda+A} h(\tau)d\tau-C}{h(\lambda+A)}, \ \ {\rm for}\ \lambda\in (0,B-A), \]
 $u\in  W^{2,p/2}(\Omega)\cap C(\bar{\Omega})$,
  $A<u<B$ a.e. in $\Omega$.

Assume further that one of the following conditions i) or ii) hold where:
\begin{description}
\item[i)]
There exists sequence $\{ \Omega_k\}_{k\in \mathbf{N}}$  of bounded  subdomains of $\Omega$ of class ${\cal C}^{0,1}$ such that $\bar{\Omega}_k\subseteq \Omega$, $\bigcup_k\Omega_k=\Omega$ and
 \[
\lim_{k\to\infty}\int_{\partial\Omega_k}\Phi_p(\nabla u(x))\cdot n(x) H_{C,A}(u(x))d\sigma (x)\in [-\infty,0],
\]
where $n(x)$ denotes unit outer normal vector to $\partial\Omega_k$, defined for $\sigma$ almost all $x\in \partial\Omega_k$;
\item[ii)] $\Omega\subseteq \mathbf{R}^n$  is Lipschitz  and
\[
\int_{\partial\Omega}\Phi_p(\nabla u(x))\cdot n(x) H_{C,A}(u(x))d\sigma (x)\in [-\infty,0],
\]
where $n(x)$ denotes unit outer normal vector to $\partial\Omega$, defined for $\sigma$ almost all $x\in \partial\Omega$.
\end{description}
Then
\[
\int_{\Omega}|\nabla u(x)|^ph(u(x))dx\leq C(n,p)\int_{\Omega} \left( \sqrt{ |\nabla^{(2)} u(x)||{\cal T}_{h,C, A}(u(x))|}\right)^{p}h(u(x))dx,
\]
where  $C(n,p)=\left( p-1+\sqrt{n-1} \right)^{\frac{p}{2}}$.
 \end{theo}

\noindent
{\bf Proof.} The proof follows directly from Theorems \ref{main1} and \ref{main2}.

Let us  substitute $u_A:= u-A:\Omega\rightarrow (0,B-A)$ a.e. instead of $u$ and $h_A(\lambda):=h(\lambda+A)$ instead of $h(\cdot)$ in Theorems \ref{main1} and \ref{main2}.
 Then $h_A(u_A)=h(u)$ and  ${\cal T}_{h_A,C}(u_A) ={\cal T}_{h,C,A}(u)$.
 In case of condition i) we  apply Theorem \ref{main1}, while in case of condition ii), we apply Theorem \ref{main2}. Verification of their assumptions is left to the reader.
 \hfill$\Box$

 Our remaining sections are devoted to the proof of the above statements and discussion.

\section{Auxiliary facts  dealing with multiplications and compositions of Sobolev functions}

The following two lemmas will be helpful for the proof of Theorem \ref{main1}.

\begin{lem}[multiplicative property of Sobolev vector field] \label{dwa}
Let $\Omega\subseteq\r^n$ be the bounded domain of class ${\cal C}^{0,1}$, $1<p<\infty$, $w_1\in W^{1,1}(\Omega,\mathbf{R}^n)\cap L^{p^{'}}(\Omega,\mathbf{R}^n)$ and $w_2\in  W^{1,p}(\Omega)\cap L^{\infty}(\Omega)$.
Then vectorfield $v:= w_1w_2$ satisfies: $v\in W^{1,1}(\Omega,\mathbf{R}^n)$ and
\begin{equation}\label{gauss}
\int_\Omega {\rm div}v (x)dx =\int_{\partial\Omega}v(x)\cdot n(x) d\sigma (x),
\end{equation}
where  $n(x)$ is outer normal vector defined $\sigma$ for almost all $x\in \partial\Omega$.
\end{lem}

\noindent
{\bf Proof.} At first we note that $v\in L^1(\Omega,\r^n)$. We  verify that $v$ is absolutely continuous along $\sigma$ almost all lines parallel to the axes and   the derivatives of the $k$-th coordinate of $v$
computed almost everywhere read as
\(
\frac{\partial v_k}{\partial x_i}(x)= \frac{\partial w_{1,k}}{\partial x_i}(x)w_2(x) + w_{1,k} \frac{\partial w_2}{\partial x_i}(x)\in L^1(\Omega).
\)
%H\"older inequality dealing with $(q,q^{'})= (1,\infty)$ and  $(q,q^{'})= (p^{'},p)$
This together with Theorem \ref{nikos} give $v\in W^{1,1}(\Omega,\r^n)$.
Trace theory (see e.g. Theorem 6.4.1 in \cite{kjf}) implies that when $w\in W^{1,1}(\Omega)$ and $\Omega$ is of class ${\cal C}^{0,1}$, then the restriction of $w$ to $\partial\Omega$ defined in the sense of trace operator belongs to $L^1(\partial\Omega, d\sigma)$.
Moreover, the formulae  \eqref{gauss}
  is the simple consequence of density of  $C^1(\bar{\Omega},\r^n)$  in $W^{1,1}(\Omega,\r^n)$, trace theorem and differentiation by parts.
This finishes the proof of the statement.
\hfill$\Box$

 \begin{lem}[compositions with Sobolev mappings]\label{jed}
 Let $n\ge 2,$ $\Omega\subseteq \r^n$ be the bounded domain of class ${\cal C}^{0,1}$, $p\ge 2$, $u\in W^{2,p/2}(\Omega)\cap L^\infty(\Omega)$.\\
 Then  $\Phi_p(\nabla u)\in W^{1,1}(\Omega,\r^n)$ where  $\Phi_p(\cdot )$ is given by \eqref{fip}.
 \end{lem}

\noindent
{\bf Proof.} We verify that for $p\ge 2$ function $\Phi_p(\lambda)$ is locally Lipschitz. According to Lemma \ref{lcomp} this implies that
$(\Phi_p(\nabla u))_l$, where $v_l$ denotes the $l$-th coordinate of vector $v$, is locally absolutely continuous for almost all lines parallel to the axes. Applying Theorem \ref{nikos}, we only have to
verify that partial derivatives of $(\Phi_p(\nabla u))_l$, computed almost everywhere, are integrable on $\Omega$.
We have for $v=\nabla u$ and almost every $x$ such that $v(x)\neq 0$
\begin{eqnarray*}
\frac{\partial}{\partial x_k}\left(\Phi_p (v) \right)_l= (p-2)|v|^{p-3}\frac{v}{|v|}\cdot \frac{\partial v}{\partial x_k}v_l + |v|^{p-2}\frac{\partial v_l}{\partial x_k}.
\end{eqnarray*}
Consequently
\begin{eqnarray}\label{osss}
|\frac{\partial}{\partial x_k}\left(\Phi_p (v) \right)|\le  (p-1)|v|^{p-2}|\nabla v|\ \ {\rm a.e.}.
\end{eqnarray}
%Let us consider three cases: a) $r>n$, b) $r=p/2$ and $p>2$, c) $r=1,p=2$.\\
In case $p>2$ we note that $v\in L^{p}(\Omega)$ by classical Gagaliardo-Nirenberg inequality \cite{ga,n1}:
\begin{eqnarray}\label{clasgn}
\|\nabla u\|_{p}\le C \| u\|_\infty^{\frac{1}{2}} \|\nabla^{(2)}u\|_{p/2}^{\frac{1}{2}} +\| u\|_\infty.
\end{eqnarray}
As $\frac{p}{p-2}=(p/2)^{'}$, we get: $|v|^{p-2}\in L^{(p/2)^{'}}(\Omega)$.
 We deduce from \eqref{osss} that for all $p\ge 2$ we have $|\frac{\partial}{\partial x_k}\left(\Phi_p (v) \right)|\in L^1(\Omega)$.
This finishes the proof of the lemma.\hfill$\Box$

%\begin{rem}\label{gnrem}\rm
% \eqref{clasgn}  with $r=p/2$ shows that when $\Omega$ is a bounded Lipschitz boundary domain and $u\in W^{2,p/2}(\Omega)\cap L^\infty (\Omega)$ then $u\in W^{1,p}(\Omega)$.
%\end{rem}

\section{Proof of Theorem \ref{main1}}

We start with the following lemma, which is the key tool for our further considerations.

\begin{lem}[precise estimate involving boundary condition]\label{goal1}
Let $n\ge 2$ and $\Omega\subseteq \mathbf{R}^n$, be the bounded domain of class ${\cal C}^{0,1}$,
$2\le p<\infty$, $0<B\le \infty$, $B$, $h,H_C, {\cal T}_{h,C}$  be as in Definition \ref{hh},  $u\in  W^{2,p/2}(\Omega)$ and there exists compact interval $[a,b]\subseteq (0,B)$
such that  $u(x)\in [a,b]$ for almost every   $x\in\Omega$.
Then we have
\begin{eqnarray*}
(I(\Omega))^{\frac{2}{p}}&:=& \left( \int_{\Omega}|\nabla u(x)|^ph(u(x))dx\right)^{\frac{2}{p}}\leq\\
&~& (p-1+\sqrt{n-1})\left( \int_{\Omega} \left( \sqrt{ |\nabla^{(2)} u(x)||{\cal T}_{h,C}(u(x))|}\right)^{p}h(u(x))dx\right)^{\frac{2}{p}}\\
&+& (I(\Omega))^{\frac{2}{p}-1}  \int_{\partial\Omega}\Phi_p(\nabla u (x))\cdot n(x) H_C(u(x)) d\sigma (x) .
\end{eqnarray*}
Moreover, the  quantities: $I(\Omega)$ and $ \int_{\partial\Omega}\Phi_p(\nabla u (x))\cdot n(x) H_C(u(x)) d\sigma (x)$ are finite.
\end{lem}

Before we prove the lemma we recall the following simple fact.

\begin{fact}(\cite{kp1})\label{faktindia}
For an arbitrary  $n\times n$ matrix $A$  and for an arbitrary
unit vector $v\in \rn$,
 we have \[ |-v^tAv+\mbox{\rm tr}\, A|\le
\sqrt{n-1}|A|.\]
\end{fact}

\noindent
We are now to prove Lemma \ref{goal1}.

\noindent
{\bf Proof of Lemma \ref{goal1}.}
We consider
\[
w_1(x):= \Phi_p(\nabla u (x)),\ \ \ w_2(x)=H_C(u(x)).
\]
Lemma \ref{jed} shows that $w_1\in W^{1,1}(\Omega,\mathbf{R}^n)$.
%Moreover, Remark \ref{gnrem} shows that $u\in W^{1,p}(\Omega)$.
We deduce that the pair $(w_1,w_2)$ obeys requirements in Lemma
 \ref{dwa}. Indeed,  by \eqref{clasgn}  and the fact that $u\in W^{2,p/2}(\Omega)\cap L^\infty (\Omega)$ we get $u\in W^{1,p}(\Omega)$.
Hence $|w_1|\le |\nabla u|^{p-1}\in L^{\frac{p}{p-1}}(\Omega)$, so that $w_1\in W^{1,1}(\Omega,\mathbf{R}^n)\cap L^{p^{'}}(\Omega,\mathbf{R}^n)$. To verify the properties of $w_2$
 we use Theorem \ref{nikos} and Lemma \ref{lcomp}.
 It allows us to deduce that $w_2\in W^{1,p}(\Omega)$ if we show  that the derivatives of $w_2$ computed almost everywhere  belong to $L^p(\Omega)$.
 For this we observe that $\nabla (H_C(u(x)))=h(u(x))\nabla u (x)$. By our assumptions $h_C(u)$ is bounded, while
 $|\nabla u|\in L^p(\Omega)$. As $H_C(u)$ is bounded as well, therefore indeed $w_2\in W^{1,p}(\Omega)\cap L^\infty(\Omega)$.
 Consequently $v:= w_1w_2\in W^{1,1}(\Omega,\r^n)$ and
\begin{eqnarray}\label{trzy}
\int_\Omega {\rm div}\left( \Phi_p(\nabla u (x))H_C(u(x))\right) dx =\int_{\partial\Omega}\Phi_p(\nabla u (x))\cdot n(x) H_C(u(x)) d\sigma (x) <\infty .
\end{eqnarray}
At the same time
\begin{eqnarray}\label{cztery}
L^1(\Omega) \ni {\rm div}\left( \Phi_p(\nabla u (x))H_C(u(x))\right)  = \Delta_p u\cdot H_C(u) + |\nabla u|^{p}h(u),
\end{eqnarray}
where $\Delta_pu = {\rm div}(\Phi_p(u))$ is the $p$-Laplacian of $u$.
Let us verify that $\Delta_pu\in L^1(\Omega)$. Direct computation gives:
\begin{eqnarray*}\label{piec}
\Delta_pu &=& (p-2)|\nabla u|^{p-2}v^t[\nabla^{(2)}u]v + |\nabla u|^{p-2}\Delta u \\
&=& |\nabla u|^{p-2} \left\{  (p-1)v^t[\nabla^{(2)}u]v + \left(-v^t[\nabla^{(2)}u]v  +\Delta u \right)\right\}\ \ {\rm a.\ e.}, \nonumber
\end{eqnarray*}
where $v=\frac{\nabla u}{|\nabla u|}$ if $\nabla u\neq 0$ and $v=0$ otherwise, $v^t$ is the transposition of $v$.
Using Fact \ref{faktindia} with $A=\nabla^{(2)}u$, we get:
\begin{eqnarray}\label{szesc}
|\Delta_pu|\le (p-1+\sqrt{n-1})|\nabla u|^{p-2}|\nabla ^{(2)}u|.
\end{eqnarray}
As $|\nabla u|^{p-2}\in L^{(\frac{p}{2})^{'}}(\Omega), |\nabla ^{(2)}u|\in L^{\frac{p}{2}}(\Omega)$, it follows that r.h.s. in \eqref{szesc} is integrable in $\Omega$, as well as $\Delta_p u\cdot H_C(u)$.
Therefore also $|\nabla u|^ph(u)$ is integrable by \eqref{cztery}.
Using \eqref{trzy}, \eqref{cztery} and \eqref{szesc} and  noting that $H_C={\cal T}_{h,C}\cdot h$ we obtain
\begin{eqnarray}
I(\Omega)&=&\int_{\Omega}|\nabla u|^ph(u)dx\nonumber\\
&\stackrel{\eqref{cztery}}{=}& \int_\Omega {\rm div}\left( \Phi_p(\nabla u (x))H_C(u(x))\right) dx -\int_\Omega \Delta_p u(x) H_C(u(x))\, dx \nonumber\\
&\stackrel{\eqref{trzy}}{\le}& \int_{\partial\Omega}\Phi_p(\nabla u (x))\cdot n(x) H_C(u(x)) d\sigma (x) + \int_\Omega |\Delta_p u(x)| |{\cal T}_{h,C}u(x)| h(u(x))\, dx\nonumber\\
&\stackrel{\eqref{szesc}}{\le}& (p-1+\sqrt{n-1})\int_{\Omega}| \nabla u|^{p-2}|\nabla ^{(2)}u| |{\cal T}_{h,C}(u)|\cdot h(u) dx\nonumber\\
& +& \int_{\partial\Omega}\Phi_p(\nabla u (x))\cdot n(x) H_C(u(x)) d\sigma (x) =:A+B.\label{plaszczka}
\end{eqnarray}
Applying H\"older inequality  we get
\begin{eqnarray*}
A&\le& (p-1+\sqrt{n-1})\left( \int_{\Omega} |\nabla u|^{p} h(u)\, dx   \right)^{1-\frac{2}{p}}\cdot \left( \int_{\Omega}
\left( \sqrt{|\nabla^{(2)} u|\cdot |{\cal T}_{h,C}(u)|}\right)^{p} h(u)\, dx    \right)^{\frac{2}{p}}\\
&=& (p-1+\sqrt{n-1})I(\Omega)^{1-\frac{2}{p}}\cdot \left( \int_{\Omega}
\left( \sqrt{|\nabla^{(2)} u|\cdot |{\cal T}_{h,C}(u)|}\right)^{p} h(u)\, dx    \right)^{\frac{2}{p}}.
\end{eqnarray*}
Now it suffices to  divide inequality \eqref{plaszczka} by $I(\Omega)^{1-\frac{2}{p}}$, but for this me must be sure that $I(\Omega)$ is finite. This is clear because we have already realized that $u\in W^{1,p}(\Omega)$ and
$h|_{[a,b]}$ is bounded.
\hfill$\Box$

\begin{rem}\label{uwauwa}\rm
In the proof of Lemma \ref{goal1} we have used the fact that $w_2(x)=H_C(u(x))$ belongs to $W^{1,p}(\Omega)$. Note that $\nabla w_2(u)=h_C(u)\nabla u$ and we have no information if  $h(u)$ is bounded if we only know that $0<u<B$.
Therefore we had to assume that $u(x)\in [a,b]$ almost everywhere. If we additionally know that $h$ is bounded near zero,  the  assumption
$0<a<u(x)<b<\infty$ a. e. can be relaxed to $0<u(x)<b<\infty$ a.e..
\end{rem}

%\begin{rem}\rm
%In case $p=2$ inequality we have $\Delta_pu=\Delta u$. Easy modification of the proof gives inequality
%\begin{eqnarray*}
%  \int_{\Omega}|\nabla u(x)|^2h(u(x))dx~&\leq&\\
% \int_{\Omega}  |\Delta u(x){\cal T}_{h,C}(u(x))|h(u(x))dx
%&+&  \int_{\partial\Omega}\Phi_p(\nabla u (x))\cdot n(x) H_C(u(x)) d\sigma (x) ,
%\end{eqnarray*}
%under the assumptions of Lemma \ref{goal1}.
%\end{rem}

\noindent
Now we prove our main statement.

\smallskip
\noindent
{\bf Proof of Theorem \ref{main1}.}
We may assume that right hand side in \eqref{niermain1} is finite and left hand side  in \eqref{niermain1} is nonzero,  as otherwise the inequality holds trivially.
 We apply Lemma \ref{goal1} to sequence of subdomains  $\Omega_k$ as in the property $(u,\Omega)$.
For this, we note that $u$ restricted to $\Omega_k$ obeys assumptions of Lemma \ref{goal1} and  we have:
\begin{eqnarray*}
(I(\Omega_k))^{\frac{2}{p}}&\le& (p-1+\sqrt{n-1})(D(\Omega_k))^{\frac{p}{2}}+ I(\Omega_k)^{\frac{2}{p}-1}B(\Omega_k),
\end{eqnarray*}
where
\begin{eqnarray*}
I(V)&:=&\int_{V}|\nabla u|^ph(u)dx \\
B(V)&:=&\int_{\partial V }\Phi_p(\nabla u (x))\cdot n (x) H_C(u(x)) d\sigma (x),\\
D(V)&:=&  \int_{V}
\left( \sqrt{|\nabla^{(2)} u|\cdot |{\cal T}_{h,C}(u)|}\right)^{p} h(u)\, dx.
\end{eqnarray*}
As $\lim_{k\to\infty} B(\Omega_k)\in (-\infty,0]$, having any $s >0$, we can assume that
\[
(I(\Omega_k))^{\frac{2}{p}}\le (p-1+\sqrt{n-1})(D(\Omega))^{\frac{2}{p}} + (I(\Omega_k))^{\frac{2}{p}-1}s.
\]
This implies that $I(\Omega_k)$ cannot converge to infinity as $k\to\infty$, so that $I(\Omega_k)\to I(\Omega)\in (0,\infty)$.
Therefore
\(
I(\Omega)^{\frac{p}{2}}\le (p-1+\sqrt{n-1})D(\Omega)^{\frac{p}{2}}+ I(\Omega)^{\frac{2}{p}-1}s,\)
which finishes the proof of the statement after we let $s$ converge to zero.
\hfill$\Box$

Following our arguments more carefully one obtains the more precise statement, which deserves the special attention.

\begin{theo}\label{goal2}
Let the assumptions of Theorem \ref{main1} be satisfied.
Then we have
\begin{eqnarray*}\label{nier11}
\left( \int_{\Omega}|\nabla u(x)|^ph(u(x))dx \right)^{\frac{2}{p}} &\leq& (p-2)\left( \int_{\Omega} \left( \sqrt{ |(Pu(x)){\cal T}_{h,C}(u(x))|}\right)^{p}h(u(x))dx\right)^{\frac{2}{p}} \\&+& \left( \int_{\Omega} \left( \sqrt{ |\Delta u(x){\cal T}_{h,C}(u(x))|}\right)^{p}h(u(x))dx\right)^{\frac{2}{p}}  ,
\end{eqnarray*}
where  $Pu$ is the nonlinear second order operator given by the formulae:
\begin{eqnarray*}
Pu(x):=\left\{
\begin{array}{ccc}
v^t\nabla^{(2)}uv, \ {\rm where}\  v=\frac{\nabla u(x)}{|\nabla u(x)|} & {\rm if}& \nabla u\neq 0\\
0 &{\rm if}& \nabla u=0.
\end{array}
\right.
\end{eqnarray*}
In particular, when $p=2$ we have
\begin{eqnarray*}
 \int_{\Omega}|\nabla u(x)|^2h(u(x))dx  &\leq&  \int_{\Omega}  |\Delta u(x)|{\cal T}_{h,C}(u(x))| h(u(x))dx .
\end{eqnarray*}
\end{theo}

\noindent
{\bf Proof.} We only have to modify the statement of Lemma \ref{goal1}, which reads as follows:
\begin{eqnarray*}
(I(\Omega))^{\frac{2}{p}}&:=& \left( \int_{\Omega}|\nabla u(x)|^ph(u(x))dx\right)^{\frac{2}{p}}\leq\\
&~& (p-2)\left( \int_{\Omega} \left( \sqrt{ |P u(x)||{\cal T}_{h,C}(u(x))|}\right)^{p}h(u(x))dx\right)^{\frac{2}{p}}\\
&+& \left( \int_{\Omega} \left( \sqrt{ |\Delta u(x)||{\cal T}_{h,C}(u(x))|}\right)^{p}h(u(x))dx\right)^{\frac{2}{p}}\\
&+& (I(\Omega))^{\frac{2}{p}-1}  \int_{\partial\Omega}\Phi_p(\nabla u (x))\cdot n(x) H_C(u(x)) d\sigma (x) .
\end{eqnarray*}
To obtain it, instead of inequality \eqref{szesc} we use:
\[
|\Delta_pu|\le (p-2)|\nabla u|^{p-2}|Pu| +|\Delta u|
\]
and substitute  it to \eqref{plaszczka}. From there the resulting inequality easily  follows.
\hfill$\Box$

\section{Proof of Theorem \ref{main2}}

%Our further analysis is restricted to the case when domain $\Omega$ is of class ${\cal C}^{1,1}$.
%Analysis within internal subdomains was needed as in general, even when we deal with the domain with regular boundary, we cannot deduce that $H_C(u)$ belongs to $W^{1,p}(\Omega)$.

In this section we prove that if $\Omega$ is Lipschitz,
we can find sequence of Lipschitz boundary subdomains $\Omega_k\subseteq \Omega$ such that $\bar{\Omega}_k\subseteq\Omega$ and
\[
\int_{\partial \Omega_k }\Phi_p(\nabla u (x))\cdot n (x) H_C(u(x)) d\sigma (x)\stackrel{k\to\infty}{\rightarrow} \int_{\partial \Omega }\Phi_p(\nabla u (x))\cdot n (x) H_C(u(x)) d\sigma (x).
\]

Our main argument is based on the following lemma. However it seems to be obvious to the specialists, it looks that its proof requires several delicate arguments. Therefore we present it in details.

\begin{lem}\label{obszary}
Assume that $\Omega$ is bounded  domain of class ${\cal C}^{0,1}$, $w\in W^{1,1}(\Omega), T\in C(\bar{\Omega})$. Then there exist sequence of  domains $\Omega_k\subseteq \Omega$ such that $\bar{\Omega}_k\subseteq\Omega$,
$\Omega_k\in {\cal C}^{0,1}$ and
\[
\int_{\partial \Omega_k }w(x)\cdot n (x) T(x) d\sigma (x)\stackrel{k\to\infty}{\rightarrow} \int_{\partial \Omega }w(x)\cdot n (x) T(x) d\sigma (x).
\]
\end{lem}

\noindent
{\bf Proof.}
The proof follows by steps.\\
{\sc Step 1.} Let us fix $t>0$, $x_0\in \mathbf{R}^n$ and define the dilation mapping ${\cal A}_{x_0,t}: \mathbf{R}^n\rightarrow \mathbf{R}^n$ by formulae
\[
{\cal A}_{x_0,t}(x) = t(x-x_0)+ x_0.
\]
We observe that:
\begin{description}
\item[1)] ${\cal A}_{x_0,t}\circ {\cal A}_{x_0,1/t}= Id$ (and so ${\cal A}_{x_0,1/t}\circ {\cal A}_{x_0,t}= Id$ i.e. the mapping ${\cal A}_{x_0,1/t}$ is inverse to ${\cal A}_{x_0,1/t}$.
\item[2)] When $V$ is any bounded domain which is starshaped with respect to $x_0\in V$ (i.e. for any $x\in V$ the segment $[x_0,x]$ is contained in $V$, see e.g. \cite{ma}, Section 1.1.8), we have
\[
\overline{{\cal A}_{x_0,t}(V)}\subseteq V\ \hbox{\rm when}\ t<1\ \hbox{\rm and}\  \overline{V}\subseteq {\cal A}_{x_0,t}(V)\ \hbox{\rm when}\ t>1.
\]
\end{description}
Indeed, properties 1) and the fact that
\begin{equation}\label{raz}
\overline{{\cal A}_{x_0,t}(V)}\subseteq V\ \hbox{\rm when}\ t<1
\end{equation}
 follow by simple verification, while the property $\overline{V}\subseteq {\cal A}_{x_0,t}(V)\ \hbox{\rm when}\ t>1$ follows
from \eqref{raz} applied to ${\cal A}_{x_0,t}(V)$ instead of $V$ as we have:
\[
\overline{V}=\overline{ {\cal A}_{x_0,1/t}\circ {\cal A}_{x_0,t}(V)}= \overline{{\cal A}_{x_0,1/t}\left( {\cal A}_{x_0,t}(V) \right)}\stackrel{\eqref{raz}}{\subseteq} {\cal A}_{x_0,t}(V).
\]

\noindent
{\sc Step 2.}
 We use the fact that every domain of class $C^{0,1}$  is finite sum of starshaped domains (see Section 1.1.9 in \cite{ma})
 and propose the formulae for the approximating family of functions.

 \smallskip
\noindent
Let $\Omega_1,\dots,\Omega_m\subseteq \Omega$ be domains  starshaped
with respect to balls $B_i=B(x_i,r_i)$ such that
  $\Omega=\cup_{i=1}^m \Omega_i$ and let us chose the appropriate smooth resolution of the unity $\{ \phi_k\}_{k=1,\dots,m}$ subordinated to $\{ \Omega_i\}$. In particular we can assume that $\phi_i$'s are smooth and compactly supported in  $\mathbf{R}^n$, moreover $\phi_i|_\Omega \equiv 0$ outside $\Omega_i$,
    $0\le\phi_i\le 1$ and $\sum_i\phi_i\equiv 1$.
  Then we define
  \[
  \Omega_t :=\cup_{k=1}^m \Omega_{k,t}, \ {\rm where}\  0<t<1,\ \Omega_{k,t}= {\cal A}_{x_{i},t}(\Omega_i).
\]
For simplicity let us denote ${\cal A}_{x_i,t}$ by $A_{i,t}$
Clearly, $\Omega_t\in {\cal C}^{0,1}$ and $\bar{\Omega}_t\subseteq \Omega$, and we have
\begin{eqnarray}\label{oznaczenia1}
I_t&:=& \int_{\partial\Omega_t} w(x)\cdot n(x) T(x)d\sigma (x) =\sum_{i=1}^m \int_{\partial\Omega_t}\phi_i (x) w(x)\cdot n(x) T(x)d\sigma (x)\nonumber\\
&=:& \sum_{i=1}^m I_{i,t},\nonumber\\
I_{i,t} &=& \int_{\partial\Omega_t\cap\Omega_i}\phi_i (x) w(x)\cdot n(x) T(x)d\sigma (x)=: \int_{\partial\Omega_t\cap\Omega_i}\phi_i (x) \Psi(x)d\sigma (x),\nonumber\\
\Psi (x) &:=& w(x)\cdot n(x)\cdot T(x)\ \hbox{\rm defined on $\partial\Omega_t$, where $t\in (0,1]$}.
\end{eqnarray}

By the very definition of $\Omega_t$, any given $x\in \partial\Omega_t\cap \Omega_i$ up to $\sigma$ measure zero comes from:
1) $A_{i,t}(y)$ for some $y\in \partial\Omega_i$,
2) $A_{j,t}(y)$ for some $y\in \partial\Omega_j$, where $j\neq i$.

More precisely, we have $\partial\Omega_t\cap \Omega_i = T_{i,t}\cup R_{i,t}\cup S_{i,t}$,  where
\begin{eqnarray*}
T_{i,t}&=& A_{i,t}\left(K_i^0\right)\ {\rm where} \ K_i^0 =   \partial\Omega_i\setminus {\Omega}\subseteq  \partial\Omega_i  \\
R_{i,t} &=& A_{i,t}\left( \partial\Omega_i\cap \Omega   \right)\cap \partial\Omega_t,\\
S_{i,t} &=& \bigcup_{j\neq i} A_{j,t}(\partial\Omega_j\setminus \partial\Omega_i)\cap \Omega_i\cap \partial\Omega_t =: \bigcup_{j\neq i} S_{i,j,t}  .
\end{eqnarray*}
In the preceding steps we will successively  analyze the integrals $\int_{K} \phi_i \Psi d\sigma$, where   $K\in \{ T_{i,t}, R_{i,t}, S_{i,t}\}$.

\noindent{\sc Step 3.} We show that $\int_{T_{i,t}} \phi_i \Psi d\sigma \rightarrow \int_{K_i^0} \phi_i \Psi d\sigma $ as $t\to 1$.

 Using change of variable formula, we obtain
\begin{eqnarray}\label{change}
\int_{T_{i,t}} \phi_i\Psi d\sigma = \int_{A_{i,t}(K_i^0) } \phi_i\Psi d\sigma = \int_{K_i^0} (\phi_i \cdot \Psi ) (A_{i,t}(x))A_{i,t}^* (d\sigma)(x),
\end{eqnarray}
where $A_{i,t}^* (d\sigma) ({ C}):= \sigma ( A_{i,t}^{-1}({ C}))= t^{-(n-1)}\sigma ({ C})$ is the pullback of the  $n-1$- dimensional Hausdorff measure under the dilation $A_{i,t}$.
Moreover, we have
\begin{eqnarray*}
\phi_i\circ A_{i,t} \rightrightarrows \phi_i \ {\rm and}\ T\circ A_{i,t} \rightrightarrows T\  {\rm as} \ t\to 1,\ \hbox{\rm uniformly on $\bar{\Omega}$},\\
\Psi (A_{i,t}(x))= w(A_{i,t}(x))\cdot n_{A_{i,t}(x)} T(A_{i,t}(x)).
\end{eqnarray*}
Applying  the telescoping argument it suffices to show that
\[
\int_{K_i^0}\phi_i (x) w(A_{i,t}(x))\cdot n_{A_{i,t}(x)} T(x) d\sigma (x) \stackrel{t\to 1}{\rightarrow} \int_{K_i^0}\phi_i (x) w(x)\cdot n_{x} T(x) d\sigma (x),
\]
where $n_{A_{i,t}(x)}$ is outer normal to $\partial\Omega_t$ at $A_{i,t}(x)$.
We  have
\begin{eqnarray*}
{\cal L}(t):= |\int_{\partial\Omega_i}\phi (x) \left( w(A_{i,t}(x))\cdot n_{A_{i,t}(x)} -  w(x)\cdot n_{x} \right)    T(x) d\sigma (x) |~~~~~~~~~~~~~~~~~~~~~~\\
~~~~~~~~~~~~~~~~~~~~~~~~~~~~~~~~~~~~~\le \| T\|_{\infty,\bar{\Omega}}\int_{\partial\Omega_i} |w(A_{i,t}(x))\cdot n_{A_{i,t}(x)} -  w(x)\cdot n_{x}|d\sigma (x).
\end{eqnarray*}
Let us consider vectors $n_x$- the outer normal to $\partial\Omega_i$ at $x$ and $n_{A_{i,t}(x)}$- the outer normal to $\partial\Omega_t$ at $A_{i,t}(x)$. Second one is the same as outer normal to $A_{i,t}(\partial\Omega_i)$ at $A_{i,t}(x)$.
It is important to note that vectors $n_x$ and $n_{A_{i,t}(x)}$ are the same. Indeed,
when the vector $w$ belongs to tangent space to $\partial\Omega_i$ at $x$, then the vector $(DA_{i,t})w=(t\cdot Id) w =tw$ belongs to tangent space at $A_{i,t}(x)$ to $A_{i,t}(\partial\Omega_i)$, because $DA_{i,t}=tId$ is the differential of the dilation $A_{i,t}$. This shows that tangent spaces $T_x\partial \Omega_i$ and $T_{A_{i,t}(x)}(A_{i,t}(\partial\Omega_i))= T_{A_{i,t}(x)}(\partial\Omega_t)$, are the same for $\sigma$ almost every $x\in K_i^0$.
As the dilation mapping $A_{i,t}$ is conformal, it does not change the angles between vectors.  Therefore  normal spaces to the respective~tangent spaces  are the same.
 This together with an easy geometric observation shows that
$n_x = n_{A_{i,t}(x)}$. We obtain
\begin{eqnarray*}
{\cal L}(t) \le
\| T\|_{\infty,\bar{\Omega}}\int_{\partial\Omega_i} |w(A_{i,t}(x)) -  w(x)|d\sigma (x)\\
\le \| T\|_{\infty,\bar{\Omega}}\int_{\partial\Omega_i} \int_t^1 |Dw(A_{i,\tau}(x)), x-x_0)| d\tau d\sigma
 \le d_i \int_t^1 \int_{\partial\Omega_i} |Dw(A_{i,\tau}(x))|d\sigma d\tau\\
\stackrel{\hbox{\rm as in \eqref{change}}}{=} d_i   \int_t^1 {\tau^{n-1}}\int_{A_{i,\tau}(\partial\Omega_i)} |Dw(x)|d\sigma d\tau
\le d_i   \int_t^1 \int_{A_{i,\tau}(\partial\Omega_i)} |Dw(x)|d\sigma d\tau \\
\le d_i C\int_{\Omega_i\setminus \Omega_{i,t}}|Dw (y)|dy \stackrel{t\to 1}{\rightarrow} 0,
\end{eqnarray*}
where $d_i= \| T\|_{\infty,\bar{\Omega}} {\rm diam}\Omega_i$.
In the last line we have used the Coarea formulae (see  \cite{fe}, Theorem 3.2.12), which can be applied as far as we prove that  $A_{i,\tau}(\partial\Omega_i)$ are level sets of the scalar Lipschitz map $\Pi_2: \bar{\Omega}\setminus \Omega_t\rightarrow \r$ constructed below. In the above notation $C$ stands for Lipschitz constant of $\Pi_2$.

It remains  to explain that $A_{i,\tau}(\partial\Omega_i)$ are level sets of some Lipschitz map. For simplicity we can assume that $x_i=0,i=1$, so that $\Omega_i=\Omega$ is starshaped with respect to some ball centered at the origin.
We define the mapping $R: \partial \Omega\times [t,1] \rightarrow \bar{\Omega}\setminus \Omega_t$ by the formulae $R(y,\tau)=\tau y$, $\tau \in [t,1]$. An easy computation shows that $R$ is Lipschitz.
%Our required argument is based on the fact that its inverse, so also its projection onto second coordinate, is Lipschitz.
We will recognize that $R^{-1}=(\Pi_1,\Pi_2)$ where $\Pi_1:\bar{\Omega}\setminus\Omega_t\rightarrow\partial\Omega$, $\Pi_2:\bar{\Omega}\setminus\Omega_t\rightarrow [t,1]$, are Lipschitz.
According to the Lemma in Section 1.1.8 in \cite{ma}, when $r=r(\omega)$ is the equation of $\partial\Omega$ in the spherical coordinates, then the function $r(\cdot)$ is Lipschitz. This implies that the mapping $\Pi_1: \bar{\Omega}\setminus \Omega_t\rightarrow \partial\Omega$, given by the formula
$\Pi_1(y):= r(\frac{y}{\| y\|})$, is also Lipschitz, because $y\mapsto \frac{1}{\| y\|}$ is positive for $y\in \bar{\Omega}\setminus \Omega_t$, and so the internal map $y\mapsto \frac{y}{\| y\|}$ is Lipschitz as well.
We construct our final mapping $\Pi_2(y):=\frac{\| y\|}{\|\Pi_1(y)\|}: \bar{\Omega}\setminus \Omega_t\rightarrow [0,t]$ and note that $\Pi_2(\cdot )$ is Lipschitz as well. Moreover, we have
$\Pi_2(y)=\tau$ if and only if  $y\in A_\tau (\partial\Omega)$. This closes our argument in the proof of Step 3.

\noindent{\sc Step 4.} We show that $\int_{R_{i,t}} \phi_i \Psi d\sigma \rightarrow 0$ as $t\to 1$.

For this, we note that $R_{i,t}\subseteq A_{i,t}\left(\partial\Omega_i\cap  \Omega\right)=: A_{i,t}(K^1_i)$.
The considerations as in the previous step and with the notation as in \eqref{oznaczenia1}, give
\begin{eqnarray*}
|\int_{R_{i,t}}\phi_i\Psi d\sigma |\le \int_{A_{i,t}(K^1_i)}\phi_i|\Psi|d\sigma \stackrel{t\to 1}{\rightarrow} \int_{K^1_i}\phi_i|\Psi|d\sigma.
\end{eqnarray*}
The latter one is zero because $\phi_i|_\Omega$ is zero outside $\Omega_i$.

\noindent{\sc Step 5.} We show that $\int_{S_{i,j,t}} \phi_i \Psi d\sigma \rightarrow 0 $ as $t\to 1$, recalling that $S_{i,j,t}:= A_{j,t}\left( \partial\Omega_j\setminus\partial\Omega_i  \right)\cap \Omega_i\cap \partial\Omega_t$.

Let us decompose
\begin{eqnarray*}
S_{i,j,t} &=&   A_{j,t}\left( \partial\Omega_j\setminus\partial\Omega_i  \right)\cap A_{j,t}\circ A_{j,t^{-1}}\left( \Omega_i\cap \partial\Omega_t\right)\\
&=& A_{j,t}\left( (\partial\Omega_j\setminus\partial\Omega_i)\cap W_{j,t}  \right) =  A_{j,t}\left(K_{i,j,t}\right)   ,\ {\rm where}\\
W_{j,t}&:=& A_{j,t^{-1}} \left( \Omega_i\cap\partial\Omega_t  \right),\\
K_{i,j,t} &:=& (\partial\Omega_j\setminus\partial\Omega_i)\cap W_{j,t}.
\end{eqnarray*}
Proceeding as in the proof of Step 4, we verify that
\begin{eqnarray*}
\int_{S_{i,j,t}}\phi_i\Psi d\sigma &=& t^{-(n-1)}\int_{K_{i,j,t}  } (\phi_i\Psi) \left( A_{j,t}(x) \right) d\sigma \\
&\le& t^{-(n-1)}\left\{ \int_{K_{i,j,t}  }|(\phi_i\Psi)( A_{j,t}(x)) - (\phi_i\Psi)(x)|d\sigma + \int_{K_{i,j,t}  }|(\phi_i\Psi)(x)|d\sigma \right\}\\
&=:& t^{-(n-1)}\left\{ X(t) + Y(t)\right\},
\end{eqnarray*}
and $X(t)\to 0$ as $t\to 1$. We are left with the estimations of second term. Observe that
\[
|Y(t)|\le \| T\|_{\infty,\bar{\Omega}}\int_{\partial\Omega_j\setminus\partial\Omega_i}|w(x)|\chi_{  W_{j,t} }(x)d\sigma (x),
\]
and we know that $w\in L^1 (\partial\Omega_j\setminus\partial\Omega_i, d\sigma)$. Now the fact that $Y(t)\to 0$ as $t\to 1$ follows from Lebesgue's Dominated Convergence Theorem after we show that
\begin{equation}\label{ostatniet}
f_t(x):= \chi_{  W_{j,t} }(x) \stackrel{t\to 1}{\rightarrow} 0,\ \hbox{\rm for every $x\in \partial\Omega_j\setminus\partial\Omega_i$}.
\end{equation}
Indeed, let $x\in \partial\Omega_j\setminus\partial\Omega_i$. If \eqref{ostatniet} was not true, we would find the sequence $t_k\to 1$ such that $f_{t_k}(x)= 1$ for every $k$.
Then we would find $y_k\in \Omega_i\cap \partial\Omega_t$ such that
$x=x_j+ t_k^{-1}( y_k-x_j)= (1-t_k^{-1})x_j + t_k^{-1}y_k.$ Consequently
\[
y_k= t_k\left( x-(1-t_k^{-1})x_j\right) \stackrel{k\to\infty}{\rightarrow} x.
\]
This is however impossible. Indeed,  by the choice of $y_k$ it would imply that $x\in \bar{\Omega}_i\cap\partial\Omega =\partial\Omega_i\cap\partial\Omega$, which contradicts the fact that
$x\not\in\partial\Omega_i$. This proves \eqref{ostatniet} and completes the arguments of Step 5, closing the whole proof of the statement.\hfill$\Box$

\noindent
{\bf Proof of Theorem \ref{main2}.} We apply Lemma \ref{obszary} to $w(x)=\Phi_p(x), T(x)=H_C(u(x))$ and use Lemma \ref{jed} and the fact that $u\in C(\bar{\Omega})$ to show that $w\in W^{1,1}(\Omega,\r^n)$ and $T\in C(\bar{\Omega})$.
\hfill$\Box$

\begin{rem} {\rm
Note that under the assumptions in Theorem \ref{main2} we cannot deduce that $H_C(u)$ belongs to $W^{1,p}(\Omega)$. This is because $\nabla H_C(u)=h(u)\nabla u$. It is clear from \eqref{clasgn} that $|\nabla u|\in L^p(\Omega)$, but $h(u)$ may be unbounded on $\Omega$. Recalling Remark  \ref{uwauwa}, even with some extra assumptions on the regularity of $\Omega$ and $u$, the proof of Theorem \ref{main2}  cannot be obtained by simplification of the proof of Theorem \ref{main1} directly with avoiding the usage of the internal subdomains $\Omega_k$.
}
\end{rem}

{\bf Acknowledgement.} This research originated when T.Ch. was  PHD student at Faculty of Mathematics, Informatics and Mechanics at the University of Warsaw. It was continued  when A.K. had research position at Institute of Mathematics of the Polish
Academy of Sciences at Warsaw. She would like to thank IM PAN for hospitality.

\end{document}